\documentclass{amsart}
\usepackage[latin2]{inputenc}
\usepackage{amsmath,amssymb,amsbsy,amsfonts,amsthm}
\usepackage{epsfig} 
\usepackage[top=3cm,bottom=2.5cm,left=2.5cm,right=2.5cm]{geometry}

\begin{document}

\numberwithin{equation}{section} \marginparwidth=2cm
\def\note#1{\marginpar{\small #1}}

\def\tens#1{\pmb{\mathsf{#1}}}
\def\vec#1{\boldsymbol{#1}}

\def\norm#1{\left|\!\left| #1 \right|\!\right|}
\def\fnorm#1{|\!| #1 |\!|}
\def\abs#1{\left| #1 \right|}
\def\ti{\text{I}}
\def\tii{\text{I\!I}}
\def\tiii{\text{I\!I\!I}}

\def\pard#1{\partial_{#1}}

\def\diver{\mathop{\mathrm{div}}\nolimits}
\def\grad{\mathop{\mathrm{grad}}\nolimits}
\def\Div{\mathop{\mathrm{Div}}\nolimits}
\def\Grad{\mathop{\mathrm{Grad}}\nolimits}

\def\tr{\mathop{\mathrm{tr}}\nolimits}
\def\cof{\mathop{\mathrm{cof}}\nolimits}
\def\det{\mathop{\mathrm{det}}\nolimits}

\def\lin{\mathop{\mathrm{span}}\nolimits}
\def\pr{\noindent \textbf{Proof: }}
\def\pp#1#2{\frac{\partial #1}{\partial #2}}
\def\dd#1#2{\frac{\d #1}{\d #2}}

\def\T{\mathcal{T}}
\def\R{\mathcal{R}}
\def\Re{\mathbb{R}}
\def\bx{\vec{x}}
\def\be{\vec{e}}
\def\bef{\vec{f}}
\def\bec{\vec{c}}
\def\bs{\vec{s}}
\def\ba{\vec{a}}
\def\bn{\vec{n}}
\def\bphi{\vec{\varphi}}
\def\btau{\vec{\tau}}
\def\bc{\vec{c}}
\def\bg{\vec{g}}

\def\bW{\tens{W}}
\def\bA{\tens{A}}
\def\bT{\tens{T}}
\def\bD{\tens{D}}
\def\bF{\tens{F}}
\def\bB{\tens{B}}
\def\bC{\tens{C}}
\def\bV{\tens{V}}
\def\bS{\tens{S}}
\def\bI{\tens{I}}
\def\bi{\vec{i}}
\def\bv{\vec{v}}
\def\bfi{\vec{\varphi}}
\def\bk{\vec{k}}
\def\b0{\vec{0}}
\def\bom{\vec{\omega}}
\def\bw{\vec{w}}
\def\p{\pi}
\def\bu{\vec{u}}
\def\bq{\vec{q}}
\def\bmu{\vec{\mu}}

\def\ID{\mathcal{I}_{\bD}}
\def\IP{\mathcal{I}_{p}}
\def\Pn{(\mathcal{P})}
\def\Pe{(\mathcal{P}^{\eta})}
\def\Pee{(\mathcal{P}^{\varepsilon, \eta})}
\def\dx{\; \mathrm{d}x}
\def\dt{\; \mathrm{d}t}

\def\Ln#1{L^{#1}_{\bn}}

\def\Wn#1{W^{1,#1}_{\bn}}

\def\Lnd#1{L^{#1}_{\bn, \diver}}

\def\Wnd#1{W^{1,#1}_{\bn, \diver}}

\def\Wndm#1{W^{-1,#1}_{\bn, \diver}}

\def\Wnm#1{W^{-1,#1}_{\bn}}

\def\Lb#1{L^{#1}(\partial \Omega)}

\def\Lnt#1{L^{#1}_{\bn, \btau}}

\def\Wnt#1{W^{1,#1}_{\bn, \btau}}

\def\Lnd#1{L^{#1}_{\bn, \btau, \diver}}

\def\Wntd#1{W^{1,#1}_{\bn, \btau, \diver}}

\def\Wntdm#1{W^{-1,#1}_{\bn,\btau, \diver}}

\def\Wntm#1{W^{-1,#1}_{\bn, \btau}}

\newtheorem{Theorem}{Theorem}[section]
\newtheorem{Example}{Example}[section]
\newtheorem{Lemma}{Lemma}[section]
\newtheorem{Proposition}{Proposition}[section]
\newtheorem{Remark}{Remark}[section]
\newtheorem{Definition}{Definition}[section]
\newtheorem{Corollary}{Corollary}[section]


\title[Regularity of parabolic systems]{$\mathcal{C}^{\alpha}$-regularity for nonlinear non-diagonal parabolic systems}\thanks{This work has been supported by the project 25-16592S financed by Czech science foundation (GA\v{C}R) and by Charles University Centre program No. UNCE/24/SCI/005. M.~Bul\'{\i}\v{c}ek is a  member of the Ne\v{c}as Center for Mathematical Modeling.}

\author{Miroslav Bul\'{\i}\v{c}ek}

\address{Mathematical Institute, Faculty of Mathematics and Physics, Charles University\\ Sokolovsk\'{a} 83,
186 75 Praha 8, Czech Republic}
\email{mbul8060@karlin.mff.cuni.cz}

\author{Jens Frehse}

\address{Institute for applied mathematics, Department of applied analysis, University of Bonn\\
Endenicher Allee 60, 53115 Bonn, Germany}
\email{erdbeere@iam.uni-bonn.de}

\begin{abstract}
In the elliptic theory for $p$-Laplacian-like problems, the H\"{o}lder continuity of solutions has been proven for problems arising as Euler--Lagrange equations of a convex potential with $p$-growth that additionally satisfies the splitting condition.

In this article, we extend these results to the parabolic setting. We investigate nonlinear parabolic systems whose structure parallels the elliptic case but incorporates time dependence. Assuming suitable space-time regularity of $F$ and natural structural conditions analogous to the stationary theory, we establish $\mathcal{C}^{\alpha}$-regularity of weak solutions in space and time whenever the growth parameter $p>d/2$. This extends the classical result for parabolic systems, which is valid only for $p>d-2$. This is the only regularity result for systems that are far from the radial (Uhlenbeck) structure.
\end{abstract}

\maketitle

\section{Introduction and statement of the result}

This paper focuses on global regularity properties  for a class of parabolic PDE's of the form
\begin{align}
\partial_t u +Lu&=f &&\textrm{in } Q, \label{2}\\
u&=0 &&\textrm{on } \Gamma, \label{bc}\\
u(0,x)&=u_0(x) &&\textrm{in } \Omega. \label{IN}
\end{align}
Here, $T>0$ is the length of time interval, $\Omega\subset \mathbb{R}^d$ is a bounded Lipschitz domain, $Q:=(0,T) \times \Omega$ denotes a parabolic cylinder and $\Gamma:=(0,T)\times \partial \Omega$. The unknown $u$ is a vector of length $N$, i.e., $u:Q\to \mathbb{R}^N$, $f:Q\to \mathbb{R}^N$ is a given right-hand side, and the elliptic operator $L$ is of the form
\begin{align}
Lu &= -\diver F_\eta (t,x,\nabla u), \label{Euler}
\end{align}
where the potential $F:Q\times \mathbb{R}^{d\times N} \to \mathbb{R}$ is a Carath\'{e}odory function satisfying certain coercivity assumptions (see below), and where we set
$$
F_{\eta}(t,x,\eta):=\frac{\partial F(t,x,\eta)}{\partial \eta} :Q \times \mathbb{R}^{d\times N} \to \mathbb{R}^{d\times N}.
$$
Here $\eta^{\alpha}_i$ corresponds to $\frac{\partial u^{\alpha}}{\partial x_i}$. We consider a class of $F$'s that for certain $p\in (1,\infty)$ provides $p$-coercivity and $p$-growth estimates, more precisely we assume that there exists $\delta_0,K\ge 0$ and $\alpha_0, \alpha_1> 0$ such that
\begin{equation}
\begin{split}
-K+\alpha_0 (\delta_0 +|\eta|^2)^{\frac{p}{2}}\le F(t,x,\eta)&\le \alpha_1(1 +|\eta|^2)^{\frac{p}{2}},\\
\left|F_{\eta}(t,x,\eta)\right| & \le \alpha_1(1+| \eta|^2)^{\frac{p-1}{2}} \label{0.0}
\end{split}
\end{equation}
for all $\eta \in \mathbb{R}^{d\times N}$ and almost all $(t,x)\in Q$. 
However, since \eqref{0.0} does not imply either convexity or quasiconvexity of $F$ with respect to $\eta$, the existence of a solution to \eqref{2}--\eqref{IN} is not guaranteed. Therefore, we also impose the $p$-convexity assumption on $F$. Namely, we assume that
\begin{align}
\alpha_0(|\eta|^2 +\delta_0)^{\frac{p-2}{2}}|\xi|^2&\le \frac{\partial^2 F(t,x,\eta)}{\partial \eta^2} \cdot (\xi \otimes \xi)\le \alpha_1(|\eta|^2 +\delta_0)^{\frac{p-2}{2}}|\xi|^2\label{5}
\end{align} 
holds true for all $\eta, \xi\in \mathbb{R}^{d\times N}$ and almost all $(t,x)\in Q$, and we complete it also by requiring certain smoothness of $F$ with respect to the spatial and time variables. To be more specific, we assume that
\begin{align}
\left|\frac{\partial F_{\eta}(t,x,\eta)}{\partial x} \right| +\left|\frac{\partial F_{\eta}(t,x,\eta)}{ \partial t} \right|  &\le \alpha_1 (\delta_0 +|\eta|^2)^{\frac{p-1}{2}}\label{5r1}\\
\left|\frac{\partial F(t,x,\eta)}{\partial x} \right|+ \left|\frac{\partial F(t,x,\eta)}{\partial t} \right| &\le \alpha_1 (1 +|\eta|^p).\label{5r2}
\end{align}
Under these assumptions, one can deduce the existence and uniqueness of a weak solution for any $p\in (1,\infty)$ provided that the initial data satisfies $u_0\in L^2(\Omega; \mathbb{R}^N)$. Moreover, the existence and uniqueness hold true even for much more complicated systems, see \cite{LaSoUr68,LaUr68}. Furthermore, one may deduce the second derivative estimates similarly as we shall do. However, further regularity results going beyond the estimates based solely on the assumptions stated above are not true. Even when 
$F$ is a convex smooth function, the global-in-time H\"{o}lder continuity of solutions is not valid, even for smooth initial data and smooth right-hand sides\footnote{For smooth data, one can expect the existence of a smooth solution only for a short time interval.}. We refer to very classical counterexamples to regularity, even in the elliptic setting, see \cite{Ne77} and \cite{SvYa02} or  \cite{DeG68}.



We summarize the available results for space-time H\"{o}lder continuity of solutions to problems under the assumptions formulated above. Under the assumptions \eqref{5}--\eqref{5r2} and sufficiently smooth data, the classical regularity results for $p\ge 2$ can be formulated for solutions $u\in L^p(0,T; W^{1,p}(\Omega; \mathbb{R}^N))$ such that it additionally satisfies  
$$
|\nabla u|^{\frac{p}{2}} \in L^{\infty}(0,T; W^{1,2}(\Omega)) \qquad \textrm{ and } \qquad \partial_t u \in L^{\infty}(0,T; L^2(\Omega; \mathbb{R}^N)).
$$
This in particular also implies (by using the classical embedding theorem)
$$
u\in L^{\infty} (0,T; W^{1,\frac{dp}{d-2}}(\Omega; \mathbb{R}^N)) \hookrightarrow L^{\infty} (0,T; \mathcal{C}^{1-\frac{d-2}{p}}(\Omega; \mathbb{R}^N)),
$$
provided $p>d-2$. Consequently, we can use the parabolic interpolation to deduce that 
\begin{equation}\label{result}
u\in \mathcal{C}^{\beta}(Q;\mathbb{R}^N) 
\end{equation}
for some $\beta>0$. 

For the better regularity results, one has to impose additional structural assumptions. Concerning the higher regularity of gradients, the only known assumption according to our knowledge is that, in addition to~\eqref{5}, we also assume
\begin{equation}
F(t,x,\eta):=F(t,x,|\eta|) \qquad \textrm{or} \qquad F(t,x,\eta):=F_0(x,|\eta|) + (o(\eta))^p, \label{Uhlen}
\end{equation}
the so-called radial (or Uhlenbeck) structure. Note that in this setting even $\mathcal{C}^{1,\alpha}$-regularity of the solution was proven by Uhlenbeck in \cite{Uh77} for the elliptic case, where the proof relies on the use of ``scalar" techniques; see also \cite{DiFr85,DiFr84} or \cite{Di93} for the parabolic setting. As mentioned above, the Uhlenbeck case \eqref{Uhlen} is the only one for which the global regularity of the solution is known. On the other hand, there is a large number of papers dealing with partial regularity of the solution under the assumption \eqref{5} only. Since we deal with everywhere $\mathcal{C}^\alpha$-regularity and not with partial regularity of the solution, we refer the interested reader to the detailed survey paper \cite{Mi06}, where the problem of partial regularity of solutions is described in great detail; see also \cite{Fu94,Gi83}. In addition, under the key assumption \eqref{Uhlen}, or some of its variants, more sophisticated estimates have recently been obtained also in the parabolic setting. The theory has been significantly extended to very sharp estimates, and we refer the interested reader to the series of papers \cite{AcMi07,Bo08,Bo14,BoDuLiSc22} and the references therein.


In the rest of the introduction, we present our results and describe the main novelty of the paper. Our primary motivation is the series of papers \cite{BuFr12,BuFrSt14,BuFrSt15}, where a new condition-the splitting condition-was introduced and for which the everywhere H\"{o}lder continuity of solutions to elliptic systems was established. We use this condition also in the present paper and extend, for a certain range of parameters $p$'s, these results to the parabolic setting as well. Hence, besides \eqref{0.0}, which is very weak, we require a kind of \(p\)-coercivity for the derivative of \(F\):
\begin{align}
\alpha_0 (\delta_0+|\eta|^2)^{\frac{p-2}{2}}|\eta|^2&\le F_\eta(t,x,\eta) \cdot \eta \le C(F(t,x,\eta)+1) \label{5.55}
\end{align}
valid for all $\eta \in \mathbb{R}^{d\times N}$  and almost all $(t,x)\in Q$.
Note that the second inequality is a consequence of \eqref{0.0}. However, in this paper we require a stronger assumption than the second inequality in \eqref{5.55}, one that more naturally reflects the $p$-setting of the problem. Thus, in the following we assume that there exists $\alpha_2>0$ and $0\le \theta<1$ such that (here $p>1$ is given by \eqref{0.0})
\begin{align}
 F_\eta(t,x,\eta) \cdot \eta \le pF(t,x,\eta) + \alpha_2(1+F(t,x,\eta))^{\theta} \label{5.5}
\end{align}
for all $\eta \in \mathbb{R}^{d\times N}$ and almost all $(t,x)\in Q$.
It seems that the presence of $p$ in \eqref{5.5} is a natural restriction, as it reflects the $p$-growth of $F$. For example,~\eqref{5.5} holds for models with the radial structure $F(t,x,\eta)\sim (1+|\eta|)^p$, as well as for more general models given, e.g., by $F(t,x,\eta)\sim (1+|\eta|)^p + (Q(\eta,\eta))^{\frac{p}{2}}$, where $Q$ denotes an arbitrary positively definite quadratic form; see also the examples below that satisfy \eqref{5.5}.

The second  and the most restrictive assumption, we assume here,  is that $F$ satisfies a splitting condition\footnote{Splitting conditions refers to the fact that $A^{\alpha \beta}$ does not depend on $i,j$. Recall, the indexes $i,j$ correspond to derivatives w.r.t. $x_i$, $x_j$ respectively.}  up to a lower order term. It means we assume that there are  symmetric matric-valued function $A: Q \times \mathbb{R}^{d\times N}\to \mathbb{R}^{N\times N}$, a symmetric matrix valued  $b: Q \to  \mathbb{R}^{d\times d}$ and  a matrix-valued function $G:Q\times \mathbb{R}^{d\times N} \to \mathbb{R}^{d\times N}$ such that
\begin{align}
F_{\eta_{i}^\alpha}(t,x,\eta)&=\sum_{\beta =1}^N \sum_{j=1}^d A^{\alpha \beta}(t,x,\eta) b_{ij}(t,x)\eta_j^\beta + G^{\alpha}_i(t,x,\eta)\label{6}
\end{align}
for all $\eta \in \mathbb{R}^{d\times N}$,  all $i=1,\ldots, d$, all $\alpha=1,\ldots, N$ and almost all  $(t,x)\in Q$. Observe that \eqref{6} is an additional structure condition which does not reflect the fact, that the $F_{\eta^{\alpha}_i}$ come from a potential. Moreover,  we require  that for all $\mu \in \mathbb{R}^{N}$, all $\eta \in \mathbb{R}^{d\times N}$, all $v\in \mathbb{R}^d$ and almost all $(t,x)\in Q$ there holds
\begin{equation}
\begin{split}
\alpha_0 (|\eta|^2 + \delta_0)^{\frac{p-2}{2}}|\mu|^2&\le\sum_{\alpha,\beta=1}^N A^{\alpha \beta}(t,x,\eta)\mu^\alpha \mu^{\beta} \le \alpha_1 (|\eta|^2 + \delta_0)^{\frac{p-2}{2}}|\mu|^2,\\
\alpha_0 |v|^2 &\le \sum_{i,j=1}^d b_{ij}(t,x) v_i v_j \le \alpha_1 |v|^2,\\
|G(\eta)|&\le \alpha_2(1+|\eta|^2)^{\frac{\theta (p-1)}{2}}, \qquad \textrm{ for some } \theta \in [0,1).\label{6.1}
\end{split}
\end{equation}

Having introduced all assumptions, we state the main result of the paper here.
\begin{Theorem}\label{main}
Let $\Omega \subset \mathbb{R}^d$ be open, and let $F$ satisfy the assumptions \eqref{0.0}, \eqref{5r2}, and \eqref{5.55}--\eqref{6.1} with $p>\max\{2,d/2\}$. Assume that $f\in W^{1,1}(0,T; L^2(\Omega; \mathbb{R}^{N}))$ and $u_0\in L^2(\Omega; \mathbb{R}^N)$. Then there exists $\alpha_{\min}>0$ such that for any weak solution $u\in L^p(0,T; W^{1,p}_0(\Omega; \mathbb{R}^N))$ to the system \eqref{2}--\eqref{Euler} and any $\delta>0$ we have
\begin{equation}\label{trub1}
u\in L^{2(p-1)}(\delta,T; \mathcal{C}^{\alpha_{\min}}(K;\mathbb{R}^N))
\end{equation}
for any compact $K \subset \Omega$. Moreover, if $u_0\in W^{1,p}(\Omega)$ then we can set $\delta=0$ in \eqref{trub1}.

If $F$ satisfies in addition the ellipticity condition \eqref{5}--\eqref{5r1}, then there exists $\beta>0$ such that the unique weak solution $u\in L^p(0,T; W^{1,p}_0(\Omega; \mathbb{R}^N))$ to the system \eqref{2}--\eqref{Euler} satisfies $u\in \mathcal{C}^{\beta}([\delta,T]\times K; \mathbb{R}^N)$ for any compact $K \subset \Omega$ and any $\delta>0$. Moreover, if $u_0 \in W^{2,2}(\Omega;\mathbb{R}^N) \cap W^{1,\infty}_0(\Omega; \mathbb{R}^N)$, then we can set $\delta=0$.

In addition, for $p\in (2,d/2]$, we have the following estimate
\begin{equation}\label{nevim3}
\sup_{t\in (0,T)} \sup_{y \in K} \int_{\Omega} \frac{|\nabla u|^p}{|x-y|^{\gamma}} \dx \le C(f,u_0,K),
\end{equation}
with
\begin{equation}\label{nevim1}
\gamma := d-p + \frac{p(2p-d)}{2(p-1)}.
\end{equation}
\end{Theorem}
We would like to emphasize here that the H\"{o}lder continuity result stated above is significantly new for $d\ge 5$. Indeed, starting from dimension $d=5$, the restriction
$$
p>\frac{d}{2}
$$
is much weaker than
$$
p>d-2,
$$
where the second restriction leads to the H\"{o}lder continuity for potentials $F$ satisfying only \eqref{5}--\eqref{5r2}. Assuming in addition the splitting-type condition, however, leads to the same result but for much lower values of $p$. We finish this part by stating a few remarks about the new splitting conditions.
We recall here a few examples \eqref{P1}--\eqref{P4} of potentials satisfying the condition \eqref{6} and refer the interested reader, e.g., to \cite{BuFr12} for more details.
 Let us define
$$
Q_\ell(\xi,\mu):=\sum_{\alpha,\beta=1}^N \sum_{i,j=1}^d a^{\alpha \beta}_\ell H_\ell^{ij}\xi_i^\alpha \mu_j^\beta,
$$
where $H_\ell$ and $a_\ell$ are symmetric positively definite constant matrices and consider for some $p\in (1,\infty)$ potentials  of the form
\begin{align}
F(\eta)&:= Q_1(\eta,\eta)|\eta|^{p-2},&&\textrm{with }H^{ij}_1:=h\delta^{ij} \textrm{ and } h>0,\label{P1}\\
F(\eta)&:=(Q_1(\eta,\eta))^{\frac{p}{2}},\label{P2}\\
F(\eta)&:= \left(Q_1(\eta,\eta)\right)^{\frac{p}{2}} + \left( Q_2(\eta,\eta)\right)^{\frac{p}{2}}, &&\textrm{with } H_1=hH_2 \textrm{ with }h\ge 0\label{P4},
\intertext{or in general of the form}
F(\eta)&:= \tilde F (Q_1(\eta,\eta),\ldots , Q_M(\eta,\eta)), &&\textrm{with } H_i=h_iH_1 \textrm{ with }h_i\ge 0 \label{P4.1}
\end{align}
for all $i=2,\ldots, M$. Furthermore, for $\tilde F$, we prescribe the following conditions: There exists $p\in (1,\infty)$ and $\{p_i\}_{i=1}^M \in \mathbb{R}$ such that $\sum_{i=1}^M p_i=p$ and
\begin{align}
\alpha_0 \prod_{i=1}^M \xi_i^{\frac{p_i}{2}} &\le \tilde F(\xi_1,\ldots, \xi_M)\le \alpha_1(1+\prod_{i=1}^M \xi_i^{\frac{p_i}{2}}),\label{til1}\\
\left|\frac{\partial \tilde F(\xi_1,\ldots, \xi_M)}{\partial \xi_k}\right| &\le \alpha_1\xi_k^{-1}\tilde F(\xi_1,\ldots,\xi_M) \textrm{ for all } k=1,\ldots, M,\label{til2}\\
\alpha_0 \prod_{i=1}^M \xi_i^{\frac{p_i}{2}}&\le 2\sum_{k=1}^M\frac{\partial \tilde F(\xi_1,\ldots, \xi_M)}{\partial \xi_k}\xi_k\le p\tilde F(\xi_1,\ldots,\xi_M). \label{til3}
\end{align}
It is  clear that we can add to each potential in \eqref{P1}--\eqref{P4.1} some lower order term (having $q$-growth with $q<p$) and that this lower order term does not change our results. It can be shown, see \cite{BuFr12}, that  all potentials given in \eqref{P1}--\eqref{P2} satisfy \eqref{0.0}, and \eqref{5.55}--\eqref{5.5}. Moreover, the potentials given by \eqref{P2}--\eqref{P4} are convex and satisfy the ellipticity condition \eqref{5}.  Concerning the splitting condition \eqref{6} on $F$, the situation is more delicate. The potential given by \eqref{P2} evidently satisfies the splitting condition. On the other hand for \eqref{P1} and \eqref{P4} we have to assume a special structure on $H$ (also given in \eqref{P1} and \eqref{P4}). To illustrate validity of the splitting condition \eqref{6}, we get for $F$ given by \eqref{P1} that
$$
F_{\eta^\alpha_i}(\eta)=\sum_{\beta=1}^N \sum_{j=1}^d 2a^{\alpha \beta}_1 H^{ij}_1\eta^\beta_j |\eta|^{p-2} + (p-2)Q_1(\eta,\eta)|\eta|^{p-4} \delta^{\alpha \beta} \delta^{ij} \eta^{\beta}_{j}.
$$
Therefore, defining
$$
A^{\alpha \beta}(t,x,\eta):= 2h a^{\alpha \beta}_1 |\eta|^{p-2} + (p-2)Q_1(\eta,\eta)|\eta|^{p-4}\delta^{\alpha \beta}, \qquad b^{ij}(x):=\delta^{ij}
$$
we see that $F$ satisfies \eqref{6} with $A^{\alpha \beta}$ and $b^{ij}$ and consequently satisfies the splitting condition. For other examples and their properties, see \cite{BuFr12}.

\section{Proof of Theorem~\ref{main}}
We provide here only the formal proof. We do not focus on the existence and uniqueness part, since it is rather standard nowadays. Furthermore, we assume that the initial condition is sufficiently regular and provide the final estimate up to the initial time $t=0$. If the initial condition $u_0$ belongs only to $L^2(\Omega;\mathbb{R}^N)$, one can adapt the method in this paper to obtain the local result on the interval $(\delta,T)$. Since such an extension is also standard, we do not present it here.

We split the proof into several parts. First, we provide the first a~priori estimates based solely on the coercivity assumptions and on the fact that the elliptic operator comes from the potential.

Second, we focus on estimates that follow from the convexity of the potential $F$. Note that these two steps are well known; therefore, we proceed quickly and formally.

Next, we deduce the weighted estimates for the gradient of the solution with the singular weight $|x-x^0|^{\gamma}$, which is the first new result of the paper. For sufficiently large $p$, we also deduce the Caccioppoli inequality, which, when combined with the weighted estimates and a proper use of the hole-filling technique, leads to regularity with respect to the spatial variable. Finally, the cross parabolic interpolation provides the main result of the paper.

\subsection{First a~priori estimates}

We recall the classical estimates for the parabolic $p$-Laplace like equations. We provide here only the formal proof. We recall that the solution is of the form $u=(u^1,\ldots, u^N)$. First, we multiply the $\alpha$-th equation in \eqref{2} by $u^{\alpha}$, take the sum over $\alpha=1,\ldots, N$, and integrate over $\Omega$. Using integration by parts and the zero boundary condition for $u$, we get
$$
\frac{\rm{d}}{\dt}\|u\|_2^2 + 2\int_{\Omega} F_{\eta} (t,x,\nabla u) \cdot \nabla u \dx = 2\int_{\Omega} f\cdot u \dx.
$$
Using the assumption \eqref{5.55}, the H\"{o}lder inequality, and the Gronwall lemma, we obtain
\begin{equation}\label{E1}
\sup_{t\in (0,T)} \|u(t)\|_2^2 + \int_Q |\nabla u|^p \dx\dt \le C\left(\|u_0\|_2,\|f\|_{L^1(L^2)}\right).
\end{equation}

Next, we repeat the same scheme with the time derivative $\partial_t u$, i.e., we multiply \eqref{2} by $\partial_t u$, integrate over $\Omega$, and, using the structure of the elliptic term, we deduce
$$ 
\int_{\Omega} |\partial_t u|^2 \dx + \frac{\rm{d}}{\dt} \int_{\Omega} F(t,x,\nabla u) \dx = \int_{\Omega} f \cdot \partial_t u + \partial_t F(t,x,\nabla u) \dx.
$$
Integration over time, together with the assumptions \eqref{0.0} and \eqref{5r2} on the data and the H\"{o}lder inequality, leads to
\begin{equation}\label{E2}
\sup_{t\in (0,T)} \|\nabla u (t)\|_p^p + \int_Q |\partial_t u|^2 \dx\dt \le C\left(\|\nabla u_0\|_p, \|f\|_{L^2(L^2)}\right).
\end{equation}
The above estimates do not require any assumption on the monotonicity of the operator $L$ - or on the convexity of the potential $F$.

\subsection{Estimates based on the convexity of $F$}
Next, we focus on the estimates that can be deduced from the ellipticity assumption \eqref{5}.  
The final regularity estimate for the time derivative is deduced in the following way. We apply the time derivative to the equation \eqref{2} and multiply the result by $\partial_t u$. Integration over $\Omega$ then leads to
\begin{equation}\label{part-blade}
\begin{split}
&\frac{\rm{d}}{\dt} \|\partial_t u\|_2^2 + 2 \int_{\Omega} \frac{\partial^2 F(t,x,\nabla u )}{\partial \eta^2}\cdot \left(\nabla \partial_t u  \otimes \nabla \partial_t u \right) \dx =  2\int_{\Omega} \partial_t f \cdot \partial_t u -\partial_t F_{\eta}(t,x,\nabla u)\cdot \nabla \partial_t u\dx\\
&\qquad \le 2\|\partial_t f\|_2\|\partial_t u\|_2 + 2\alpha_1 \int_{\Omega} (\delta_0 + |\nabla u|^2)^{\frac{p-1}{2}} |\nabla \partial_t u|\dx\\
&\qquad \le 2\|\partial_t f\|_2\|\partial_t u\|_2 + \int_{\Omega} \frac{\partial^2 F(t,x,\nabla u )}{\partial \eta^2}\cdot \left(\nabla \partial_t u  \otimes \nabla \partial_t u \right) \dx + C\int_{\Omega} (\delta_0 + |\nabla u|^2)^{\frac{p}{2}}\dx,
\end{split}
\end{equation}
where for the first inequality we used the H\"{o}lder inequality and the assumption~\eqref{5r1}, while for the second inequality, we used the Young inequality combined with \eqref{5}. Next, we can absorb the second term on the right-hand side by the corresponding term on the left-hand side. To apply the Gronwall lemma and integrate the resulting inequality with respect to time in order to obtain the final bound, we need to control $\partial_t u(0)$. The initial value of the time derivative can be estimated from the equation \eqref{2} as
$$
\|\partial_t u(0)\|_2 = \left\|f(0)+ \diver F_{\eta}(\nabla u_0)\right\|_2.
$$
Consequently, inserting this into \eqref{part-blade}, and using the assumption \eqref{5} together with the already obtained a~priori estimates \eqref{E2}, we arrive at the estimate
\begin{equation}\label{E3}
\sup_{t\in (0,T)} \|\partial_t u\|_2^2 + \int_Q (1+|\nabla u|)^{p-2}|\nabla \partial_t u|^2 \,\mathrm{d}x \,\mathrm{d}t \le C\left(\|\nabla u_0\|_{\infty}, \|\nabla^2 u_0\|_2, \|f\|_{W^{1,2}(L^2)}\right).
\end{equation}

\bigskip

The final classical regularity estimate is for the second spatial derivative. Since we do not consider any smoothness of the boundary, we provide here only the interior estimates. We multiply \eqref{2} by $-\diver (\nabla u \zeta^2)$, where $\zeta$ is a smooth nonnegative function compactly supported in $\Omega$. By integration over $\Omega$ and using integration by parts, we deduce
$$
\int_{\Omega} \frac{\partial^2F(t,x,\nabla u )}{\partial \eta^n}\cdot \left(\nabla^2 u \otimes \nabla (\nabla u \zeta^2)\right) \dx  = \int_{\Omega} (\partial_t u - f)\cdot \diver (\nabla u \zeta^2) -\frac{\partial F_{\eta}(t,x,\nabla u)}{\partial x} \cdot \nabla (\nabla u \zeta^2)\dx. 
$$
The assumption \eqref{5} and \eqref{5r1}, together with the H\"{o}lder and Young inequalities, and using also the fact that $p\ge 2$, lead to the estimate
$$
\begin{aligned}
&\alpha_0 \int_{\Omega} (1+|\nabla u|^2)^{\frac{p-2}{2}} |\nabla^2 u|^2 \zeta^2 \dx \\
&\quad \le C(\zeta) \int_{\Omega} (1+ |\nabla u|^2)^{\frac{p-1}{2}} |\nabla^2 u|\zeta  + (|f|+|\partial_t u|)|\nabla u|+(1+|\nabla u|^2)^{\frac{p}{2}}  + (|f|+|\partial_t u|) |\nabla^2 u|\zeta \dx\\
& \quad \le \frac{\alpha_0}{2} \int_{\Omega} (1+ |\nabla u|^2)^{\frac{p-2}{2}} |\nabla^2 u|^2\zeta^2 \dx + C\left(1+\|\nabla u\|_p^p + \|f\|^2_2 +\|\partial_t u\|_2^2\right).
\end{aligned}
$$
Moving the first term on the right-hand side to the left, we obtain 
$$
\begin{aligned}
&\int_{\Omega} (1+|\nabla u|^2)^{\frac{p-2}{2}} |\nabla^2 u|^2 \zeta^2 \dx  \le  C(\zeta) \left(1+\|f\|_2^2 + \|\nabla u\|^p_p + \|\partial_t u\|_2^2\right).
\end{aligned}
$$
Thus, using the a~priori estimate \eqref{E3} and the assumption on the data, we deduce that
\begin{equation}\label{E5}
\sup_{t\in (0,T)} \int_{\Omega} (1+|\nabla u|^2)^{\frac{p-2}{2}} |\nabla^2 u|^2 \zeta^2 \dx \le C\left(\zeta,\|\nabla u_0\|_{\infty}, \|\nabla^2 u_0\|_2, \|f\|_{W^{1,2}(L^2)}\right).
\end{equation}
Note that the above estimate together with \eqref{E3} can be summarized as 
\begin{equation}\label{Esum}
\begin{aligned}
\sup_{t\in (0,T)} \int_{\Omega} |\partial_t u|^2+\left|\nabla (1+|\nabla u|^2\zeta^2)^{\frac{p}{4}}\right|^2 \dx&\le C,
\end{aligned}
\end{equation}
where $C$ depends on the data and on the choice of the function $\zeta$. 

\subsection{Estimates based on the structure of $F$}
In this section, we use the splitting condition to deduce weighted estimates for $|F|$ and consider the case $1 \le p \le d$. We again employ a nonnegative, compactly supported function~$\zeta$ to localize the argument and to remain inside $\Omega$. Throughout the section we use the Einstein summation convention for repeated indices, and, to simplify the notation, we write $D_j$ instead of~$\frac{\partial}{\partial x_j}$.  

For arbitrary $\gamma\in [0,d-p]$ and $y\in \mathbb{R}^d$, we multiply the $\nu$-th equation in \eqref{2} by $D_j u^{\nu} (x_j - y_j) |x - y|^{-\gamma} \zeta$, integrate over $\Omega$ and use integration by parts and sum over $\nu=1,\ldots, N$ to get  (here $x^0:(0,T)\to \Omega$ is arbitrary).
\begin{equation} \label{basis}
\begin{split}
&\int_{\Omega} \partial_t u^{\nu} \frac{D_j u^{\nu}  (x_j-y_j)}{ |x-y|^{\gamma}} \zeta \dx + \int_{\Omega} F_{\eta^{\nu}_i}(t,x,\nabla u) D_i\left(\frac{D_j u^{\nu} (x_j-y_j)}{ |x-y|^{\gamma}} \zeta \right) \dx \\
&\qquad = \int_{\Omega} f^{\nu} \frac{D_j u^{\nu}  (x_j-y_j(t))}{ |x-y|^{\gamma}} \zeta \dx.
\end{split}
\end{equation}
In what follows we rewrite and estimate all terms appearing in \eqref{basis}.

First, we evaluate the second term on the left hand side of \eqref{2} as follows (we use integration by parts)
$$
\begin{aligned}
&-\int_{\Omega} F_{\eta^{\nu}_i}(t,x,\nabla u) D_i\left(\frac{D_j u^{\nu}  (x_j-y_j)}{ |x-y|^{\gamma}} \zeta\right) \dx=-\int_{\Omega} \frac{D_j F(t,x,\nabla u) (x_j-y_j)\zeta}{ |x-y|^{\gamma}} \dx\\
&\quad -\int_{\Omega} \frac{F_{\eta}(t,x,\nabla u) \cdot \nabla u }{ |x-y|^{\gamma}} \zeta \dx +\gamma\int_{\Omega} F_{\eta^{\nu}_i}(t,x,\nabla u) \left(\frac{D_j u^{\nu} (x_j-y_j)(x_i-y_i)}{ |x-y|^{\gamma+2}} \zeta\right) \dx\\
&\quad  -\int_{\Omega} F_{\eta^{\nu}_i}(t,x,\nabla u) \left(\frac{D_j u^{\nu}  (x_j-y_j)}{ |x-y|^{\gamma}} D_i\zeta\right) \dx+ \int_{\Omega} \frac{\partial F(t,x,\nabla u)}{\partial x_j} \frac{(x_j-y_j)\zeta}{ |x-y|^{\gamma}} \dx\\
&=(d-p-\gamma)\int_{\Omega}  \frac{F(t,x,\nabla u)\zeta}{|x-y|^{\gamma}}  \dx+\gamma\int_{\Omega} F_{\eta^{\nu}_i}(t,x,\nabla u) \left(\frac{D_j u^{\nu} (x_j-y_j)(x_i-y_i)}{ |x-y|^{\gamma+2}} \zeta\right) \dx\\
&\quad-\int_{\Omega} \frac{F_{\eta}(t,x, \nabla u)\cdot \nabla u - pF(t,x,\nabla u) }{ |x-y|^{\gamma}} \zeta \dx -\int_{\Omega} F_{\eta^{\nu}_i}(t,x,\nabla u) \left(\frac{D_j u^{\nu}  (x_j-y_j)}{ |x-y|^{\gamma}} D_i\zeta\right) \dx\\
&\quad +\int_{\Omega}  F(t,x,\nabla u) \frac{(x_j-y_j)}{ |x-y|^{\gamma}} D_j\zeta \dx+ \int_{\Omega} \frac{\partial F(t,x,\nabla u)}{\partial x_j} \frac{(x_j-y_j)\zeta}{ |x-y|^{\gamma}} \dx
\end{aligned}
$$
We use this identity in the equation~\eqref{basis} to deduce
\begin{equation} \label{basis2}
\begin{split}
&(d-p-\gamma)\int_{\Omega}  \frac{F(t,x,\nabla u)\zeta}{|x-y|^{\gamma}}  \dx+\gamma\int_{\Omega} F_{\eta^{\nu}_i}(t,x,\nabla u) \left(\frac{D_j u^{\nu} (x_j-y_j)(x_i-y_i)}{ |x-y|^{\gamma+2}} \zeta\right) \dx \\
&\quad = -\int_{\Omega} f^{\nu} \frac{D_j u^{\nu}  (x_j-y_j)}{ |x-y|^{\gamma}} \zeta \dx+\int_{\Omega} \partial_t u^{\nu} \frac{D_j u^{\nu}  (x_j-y_j)}{ |x-y|^{\gamma}} \zeta \dx\\
&\qquad +\int_{\Omega} \frac{F_{\eta}(t,x,\nabla u)\cdot \nabla u - pF(t,x,\nabla u) }{ |x-y|^{\gamma}} \zeta \dx +\int_{\Omega} F_{\eta^{\nu}_i}(t,x,\nabla u) \left(\frac{D_j u^{\nu}  (x_j-y_j)}{ |x-y|^{\gamma}} D_i\zeta\right) \dx\\
&\qquad -\int_{\Omega}  F(t,x,\nabla u) \frac{(x_j-y_j)}{ |x-y|^{\gamma}} D_j\zeta\dx - \int_{\Omega} \frac{\partial F(t,x,\nabla u)}{\partial x_j} \frac{(x_j-y_j)\zeta}{ |x-y|^{\gamma}} \dx.
\end{split}
\end{equation}
To estimate the terms on the left-hand side from below, we use the assumptions \eqref{0.0} and \eqref{6.1}. For the third term on the right-hand side, we employ the assumption \eqref{5.5}, and for the remaining terms on the right-hand side we use \eqref{0.0} and \eqref{5r2}. The above inequality then reduces to the following.
\begin{equation} \label{basis3}
\begin{split}
(d-p-\gamma)\int_{\Omega}  \frac{|\nabla u|^p \zeta}{|x-y|^{\gamma}}  \dx+&\gamma\int_{\Omega} \frac{(\delta_0+|\nabla u|^2)^{\frac{p-2}{2}} |\nabla u \cdot (x-y)|^2}{ |x-y|^{\gamma+2}} \zeta \dx \\
&\quad \le C\int_{\Omega}  \frac{(|f|+|\partial_t u|) |\nabla u \cdot (x-y)|\zeta}{ |x-y|^{\gamma}} \dx + C\int_{\Omega}\frac{(1+|\nabla u|^p)^{\theta}}{ |x-y|^{\gamma}} \zeta \dx \\
&\qquad +C\int_{\Omega} \frac{(1+|\nabla u|^p)(\zeta + |\nabla \zeta|)}{ |x-y|^{\gamma-1}}  \dx.
\end{split}
\end{equation}
In what follows, we focus on the weighted estimates based on \eqref{basis3}. 
We start with estimating the first  term on the right-hand side. We use the H\"{o}lder inequality and the fact that $p\ge 2$ as follows.
$$
\begin{aligned}
&C\int_{\Omega}  \frac{(|f|+|\partial_t u|) |\nabla u \cdot (x-y)|\zeta}{ |x-y|^{\gamma}} \dx\\
&\qquad \le C\int_{\Omega}  \frac{(|f|+|\partial_t u|)\zeta^{\frac{1}{p'}}}{|x-y|^{\frac{\gamma}{p'}-1}} \left(\frac{(\delta_0 +|\nabla u|^2)^{\frac{p-2}{2}} |\nabla u \cdot (x-y)|^2\zeta}{ |x-y|^{\gamma+2}}\right)^{\frac{1}{p}} \dx\\
&\qquad \le C(\|f\|_2 + \|\partial_t u\|_2) \left(\int_{\Omega} \frac{(\delta_0 +|\nabla u|^2)^{\frac{p-2}{2}} |\nabla u \cdot (x-y)|^2\zeta}{ |x-y|^{\gamma+2}} \dx \right)^{\frac{1}{p}} \left(\int_{\Omega}  \frac{\zeta \dx}{|x-y|^{\frac{\gamma -p'}{p'}\frac{2p}{p-2}}}  \right)^{\frac{p-2}{2p}}\\
&\qquad \le C(\|f\|_2 + \|\partial_t u\|_2)^{p'} \left(\int_{\Omega}  \frac{\zeta \dx}{|x-y|^{\frac{\gamma -p'}{p'}\frac{2p}{p-2}}}  \right)^{\frac{p'(p-2)}{2p}} + \frac{\gamma \alpha_0}{2} \int_{\Omega} \frac{(\delta_0 +|\nabla u|^2)^{\frac{p-2}{2}} |\nabla u \cdot (x-y)|^2\zeta}{ |x-y|^{\gamma+2}} \dx  .
\end{aligned}
$$
In order for the right-hand side to be finite, it is evident that we must require (recall also that $p\ge 2$ here)
\begin{equation}\label{gamma1}
\frac{\gamma -p'}{p'}\frac{2p}{p-2}< d \quad \iff \quad \gamma < \frac{d(p-2)+2p}{2(p-1)} = d-p + \frac{p(2p-d)}{2(p-1)} 
\end{equation}
or that $\gamma \le 2$ for $p=2$. Simultaneously, we need that $\gamma\le d-p$. Comparing it with \eqref{gamma1}, we see that the final restriction, we get is
\begin{equation}
\gamma \left\{ \begin{aligned}&\le d-p &&\textrm{ for } d\le 4,\\
&\le d-p &&\textrm{ for } d>4 \textrm{ and } p> \frac{d}{2}\\
&< d-p + \frac{p(2p-d)}{2(p-1)}   &&\textrm{ otherwise.}
\end{aligned} \right.
\label{gamma2}
\end{equation}

We substitute the above inequality into \eqref{basis3} and use the Young inequality to absorb the corresponding term on the right-hand side. Hence, we deduce 
\begin{equation} \label{basis4}
\begin{split}
&(d-p-\gamma)\int_{\Omega}  \frac{|\nabla u|^p \zeta}{|x-y|^{\gamma}}  \dx+\int_{\Omega} \frac{(\delta_0+|\nabla u|^2)^{\frac{p-2}{2}} |\nabla u \cdot (x-y)|^2}{ |x-y|^{\gamma+2}} \zeta \dx \\
&\quad \le C(\|f\|_2 + \|\partial_t u\|_2)^{p'} \left(\int_{\Omega}  |x-y|^{-\frac{\gamma -p'}{p'}\frac{2p}{p-2}}\zeta \dx  \right)^{\frac{p'(p-2)}{2p}}  \\
&\qquad \ \qquad  +C\int_{\Omega} \frac{(1+|\nabla u|^p)(\zeta+ |\nabla \zeta|)}{|x-y|^{\gamma-1}} \dx + C\int_{\Omega} \frac{(1+|\nabla u|^p)^{\theta}\zeta}{ |x-y|^{\gamma}} \dx,
\end{split}
\end{equation}
which is the starting point of further estimates. First, we obtain the bound whenever $\gamma < d-p$ and satisfies also \eqref{gamma2}. Indeed, in this case we can use the Young inequality to get the bound on the last term in \eqref{basis4} as follows. 
$$
\begin{aligned}
C\int_{\Omega} \frac{(1+|\nabla u|^p)^{\theta}\zeta}{ |x-y|^{\gamma}} \dx &\le  \frac{d-p-\gamma}{4} \int_{\Omega} \frac{|\nabla u|^p\zeta}{ |x-y|^{\gamma}}\dx + C(\gamma,p,d) \int_{\Omega} \frac{\zeta}{ |x-y|^{\gamma}}\dx. 
\end{aligned}
$$
Similarly, we can estimate the penultimate term on the right-hand side of \eqref{basis4} as follows. We use the fact that
$$
C\int_{\Omega} \frac{(1+|\nabla u|^p)\zeta }{|x-y|^{\gamma-1}} \dx \le \frac{d-p-\gamma}{4} \int_{\Omega} \frac{|\nabla u|^p\zeta}{ |x-y|^{\gamma}}\dx + C(\gamma,p,d) \int_{\Omega} |\nabla u|^p \dx.
$$
Inserting these two inequalities into \eqref{basis4} gives
\begin{equation} \label{basis5}
\begin{split}
(d-p-\gamma)\int_{\Omega}  \frac{|\nabla u|^p \zeta}{|x-y|^{\gamma}}  \dx &\le  C(\|f\|_2 + \|\partial_t u\|_2)^{p'} \left(\int_{\Omega}  |x-y|^{-\frac{\gamma -p'}{p'}\frac{2p}{p-2}}\zeta \dx  \right)^{\frac{p'(p-2)}{2p}} \\
&\qquad + C\int_{\Omega} \frac{\zeta}{ |x-y|^{\gamma}}+\frac{(1+|\nabla u|^p) (\zeta+|\nabla \zeta|)}{ |x-y|^{\gamma-1}} \dx.
\end{split}
\end{equation}
Hence, we assume that $K\subset \Omega' \subset   \overline{\Omega'}\subset \Omega$ are arbitrary such that $K$ is closed and $\Omega'$ is open and that $\zeta\equiv 1$ in $\Omega'$. Furthermore, assume that $y\in K$. Then the above inequality can be reduced once more as (noting that $\gamma<d-p$)
\begin{equation*}
\begin{split}
(d-p-\gamma)\int_{\Omega}  \frac{|\nabla u|^p \zeta}{|x-y|^{\gamma}}  \dx &\le  C(K)(1+\|f\|_2 + \|\partial_t u\|_2)^{p'} + C(K)\int_{\Omega}|\nabla u|^p \dx.
\end{split}
\end{equation*}
Since the right-hand side is independent of $y$, we can therefore deduce that 
\begin{equation} \label{basis6}
\begin{split}
(d-p-\gamma)\sup_{y\in K} \int_{\Omega}  \frac{|\nabla u|^p \zeta}{|x-y|^{\gamma}}  \dx &\le  C(K)(1+\|f\|_2 + \|\partial_t u\|_2)^{p'} + C(K)\int_{\Omega}|\nabla u|^p \dx.
\end{split}
\end{equation}
\subsubsection{No ellipticity condition} In the case where we do not assume \eqref{5}, we have only the estimates \eqref{E1} and \eqref{E2} at our disposal. Therefore, applying the $2/p'$ power to \eqref{basis6}, integrating over $t\in (0,T)$, and using these estimates leads to
\begin{equation} \label{basis7}
\begin{split}
\int_0^T\left(\sup_{y\in K} \int_{K}  \frac{|\nabla u|^p}{|x-y|^{\gamma}}  \dx\right)^{\frac{2}{p'}} \dt  &\le  C(K, \gamma, u_0, f)
\end{split}
\end{equation}
valid for all $\gamma$ satisfying \eqref{gamma2} and fulfilling also $\gamma < d-p$. 

\subsubsection{Estimate with ellipticity condition}
Here, we can use the a~priori estimate \eqref{E3}, which combined with \eqref{basis6} leads to
\begin{equation} \label{ECg}
\begin{split}
\sup_{t\in (0,T)} \sup_{y\in K} \int_{\Omega}  \frac{|\nabla u|^p \eta}{|x-y|^{\gamma}}  \dx &\le  C(K, u_0,f,\gamma)
\end{split}
\end{equation}
for all $\gamma < d-p$ fulfilling also \eqref{gamma2}. 

\subsection{Estimates for critical $\gamma=d-p$}
Next, if $d$ and $p$ satisfy the restrictions in \eqref{gamma2}$_1$ or \eqref{gamma2}$_2$, which means that we certainly have $p\in \left(\frac{d}{2},d\right]$, we can improve our result as follows. For any $y\in K$ and any $R$ with $B_{2R}\subset \Omega$, we choose $\zeta$ such that $\zeta\in \mathcal{C}^1_0 ( B_{2R}(y))$, $\zeta\equiv 1$ on $B_{R}(y)$, and $|\nabla \zeta|\le C/R$. Moreover, we set $\gamma:=d-p$. In what follows, we abbreviate $A_R(y):= B_{2R}(y)\setminus B_{R}(y)$. It then follows from \eqref{basis4} that
\begin{equation} \label{B10}
\begin{split}
&\int_{B_R(y)} \frac{(\delta_0+|\nabla u|^2)^{\frac{p-2}{2}} |\nabla u \cdot (x-y)|^2}{ |x-y|^{d-p+2}} \dx \\
&\quad \le C(\|f\|_2 + \|\partial_t u\|_2)^{p'} \left(\int_{B_{2R}(y)}  |x-y|^{-\frac{(d-p) -p'}{p'}\frac{2p}{p-2}} \dx  \right)^{\frac{p'(p-2)}{2p}}  \\
&\qquad \ \qquad  +C\int_{A_R(y)} \frac{(1+|\nabla u|^p)}{R^{d-p}} \dx + C\int_{B_{2R}(y)} \frac{(1+|\nabla u|^p)^{\theta}}{ |x-y|^{d-p}} \dx+ C\int_{B_{2R}(y)} \frac{1+|\nabla u|^p}{ |x-y|^{d-p-1}} \zeta\dx
\end{split}
\end{equation}
We focus on the individual terms on the right-hand side. The easiest one is the integral in the first term, for which we have
\begin{equation}\label{koule}
\left(\int_{B_{2R}(y)}  |x-y|^{-\frac{(d-p) -p'}{p'}\frac{2p}{p-2}} \, \mathrm{d}x  \right)^{\frac{p'(p-2)}{2p}} \le C R^{\alpha_0},
\end{equation}
where $\alpha_0$ is given by
\begin{equation}\label{A0}
\alpha_0:= \frac{p(2p-d)}{2(p-1)}>0,
\end{equation}
and the positivity follows from the fact that $p\in \left(\frac{d}{2},\, d\right)$. Next, to estimate the last term in~\eqref{B10}, we use~\eqref{basis5} with~$\gamma:=d-p-1$. In this way, we obtain
\begin{equation*} 
\begin{split}
\int_{\Omega}  \frac{|\nabla u|^p \zeta}{|x-y|^{d-p-1}}  \dx &\le  C(\|f\|_2 + \|\partial_t u\|_2)^{p'} \left(\int_{\Omega}  |x-y|^{-\frac{d-p-1 -p'}{p'}\frac{2p}{p-2}}\zeta \dx  \right)^{\frac{p'(p-2)}{2p}} \\
&\qquad + C\int_{\Omega} \frac{\zeta}{ |x-y|^{d-p-1}}+\frac{(1+|\nabla u|^p) (\zeta+|\nabla \zeta|)}{ |x-y|^{d-p-2}} \dx\\
&\le  C(\|f\|_2 + \|\partial_t u\|_2)^{p'} R^{\alpha_0+1}+CR^{p+1}+ C\int_{A_R(y)} \frac{|\nabla u|^p}{R^{d-p-1}}\dx \\
&\qquad + CR\int_{\Omega} \frac{(1+|\nabla u|^p)\zeta}{ |x-y|^{d-p-1}} \dx.
\end{split}
\end{equation*}
Next, we fix $R_0$ such that $C R_0 \le \frac{1}{2}$. Then, assuming that $R \le R_0$, we can absorb the last term on the right-hand side into the left-hand side to obtain
\begin{equation} \label{basis5.1}
\begin{split}
\int_{\Omega}  \frac{|\nabla u|^p \zeta}{|x-y|^{d-p-1}}  \dx &\le  C(1+\|f\|_2 + \|\partial_t u\|_2)^{p'} R^{\alpha_0+1}+ C\int_{A_R(y)} \frac{|\nabla u|^p}{R^{d-p}}\dx.
\end{split}
\end{equation}

The last term we need to estimate in \eqref{B10} is the one with the $\theta$ power. We use the H\"{o}lder inequality and the estimate \eqref{basis6} with $\gamma:= d-p- \frac{p(1-\theta)}{2\theta}<d-p$, together with the a~priori estimate \eqref{E2}, to deduce
\begin{equation}\label{B1}
\begin{split}
\int_{B_{2R}(y)} \frac{(1+|\nabla u|^p)^{\theta}}{ |x-y|^{d-p}} \dx &\le  \int_{B_{2R}(y)} \left(\frac{(1+|\nabla u|^p)}{ |x-y|^{d-p-\frac{p(1-\theta)}{2\theta}}} \right)^{\theta} \frac{1} {|x-y|^{(1-\theta)(d-p)+\frac{p(1-\theta)}{2}}} \dx\\
&\le  C\left(\int_{B_{2R}(y)} \frac{(1+|\nabla u|^p)}{ |x-y|^{d-p -\frac{p(1-\theta)}{2\theta}}} \dx \right)^{\theta} R^{\frac{p(1-\theta)}{2}}\\
&\le CR^{\frac{p(1-\theta)}{2}}\left(1+\|f\|_2 + \|\partial_t u\|_2\right)^{p'}.
\end{split}
\end{equation}
Finally, substituting \eqref{koule}, \eqref{basis5.1}, and \eqref{B1} into \eqref{B10}, we obtain
\begin{equation} \label{B2}
\begin{split}
&\int_{B_R(y)} \frac{(\delta_0+|\nabla u|^2)^{\frac{p-2}{2}} |\nabla u \cdot (x-y)|^2}{ |x-y|^{d-p+2}} \dx \\
&\quad \le C(1+\|f\|_2 + \|\partial_t u\|_2)^{p'}R^{\alpha_1} + C\int_{A_R(y)} \frac{(1+|\nabla u|^p)}{R^{d-p}} \dx,
\end{split}
\end{equation}
where $\alpha_1:=\min\{\alpha_0, \frac{p(1-\theta)}{2} \}$.

\subsection{Caccioppoli inequality, case $p\ge 2$} 
In this section, we use the assumption on the splitting condition, see \eqref{6} and \eqref{6.1}. 
For simplicity, we consider here only the case when $b_{ij}=\delta_{ij}$. 
We keep the notation from the previous part, and in addition, we assume a specific form of $\zeta$, namely that $\zeta(x) = \xi(|x-y|)$ so that consequently $\nabla \zeta = \frac{(x-y)}{|x-y|} \xi'(|x-y|)$ with $|\xi'|\le R$ and $\xi'$ supported on $A_{R}(y)$. 
Moreover, we denote
$$
\overline{u}_R:= \frac{1}{|A_R(y)|}\int_{A_R(y)} u \dx,
$$
the mean value of $u$ over the annulus $A_R(y)$. 
Finally, we multiply \eqref{2} by $(u-\overline{u}_R)\eta$ and integrate over $\Omega$ to get
\begin{equation}\label{C1}
\int_{\Omega} F_{\eta} (t,x,\nabla u) \cdot \nabla u \, \zeta  \dx = \int_{\Omega} (f-\partial_t u)\cdot (u-\overline{u}_R)\zeta \dx - \int_{\Omega}  F_{\eta} (t,x,\nabla u) \cdot  ((u-\overline{u}_R)\otimes  \nabla \zeta)\dx.
\end{equation}
The left-hand side is estimated with the help of \eqref{5.55}, the first term on the right-hand side by using the H\"{o}lder inequality, and finally the last term by using the splitting condition \eqref{6}, the assumption \eqref{6.1}, and the Cauchy--Schwarz inequality. The resulting inequality from \eqref{C1} is of the form, and we also use the structure of $\zeta$
\begin{equation}\label{C2}
\begin{split}
\int_{B_{R}(y)} (\delta_0+|\nabla u|^2)^{\frac{p-2}{2}} |\nabla u|^2 \dx &\le  CR^{\frac{d(p-2)}{2p}}(\|f\|_2 + \|\partial_t u\|_2) \left(\int_{B_{2R}(y)}|u-\overline{u}_R|^p \dx \right)^{\frac{1}{p}} \\
&+ CR^{-2}\int_{A_R(y)} (\delta_0+|\nabla u|^2)^{\frac{p-2}{2}} |\nabla u \cdot (x-y)||u-\overline{u}_R|\dx \\
& + CR^{-1}\int_{A_R(y)} (\delta_0+ |\nabla u|)^{\theta(p-1)} |u-\overline{u}_R| \dx.
\end{split}
\end{equation}
Using the Poincar\'{e} inequality and the H\"{o}lder inequality, we can estimate the terms appearing on the right-hand side as follows. For the first term, we have 
$$
\begin{aligned}
CR^{\frac{d(p-2)}{2p}}(\|f\|_2 + \|\partial_t u\|_2) \left(\int_{B_{2R}(y)}|u-\overline{u}_R|^p \dx \right)^{\frac{1}{p}}  \le CR^{d-p +\frac{2p-d}{2}}(\|f\|_2 + \|\partial_t u\|_2) \left(\int_{B_{2R}(y)}\frac{|\nabla u|^p}{R^{d-p}} \dx \right)^{\frac{1}{p}}.
\end{aligned}
$$
The second term is estimated via the H\"{o}lder and the Poincar\'{e} inequalities as
$$
\begin{aligned}
&R^{-2}\int_{A_R(y)} (\delta_0+|\nabla u|^2)^{\frac{p-2}{2}} |\nabla u \cdot (x-y)||u-\overline{u}_R|\dx\le \\
&\le\quad R^{-2}\left(\int_{A_R(y)} (\delta_0+|\nabla u|^2)^{\frac{p-2}{2}} |\nabla u \cdot (x-y)|^2 \dx \right)^{\frac12}  \cdot \\
&\qquad \ \qquad \cdot   \left(\int_{B_{2R}(y)} |u-\overline{u}_R|^p\dx \right)^{\frac{1}{p}} \left(  \int_{B_{2R}(y)} (\delta_0+|\nabla u|^2)^{\frac{p}{2}}\dx \right)^{\frac{p-2}{2p}}\\
&\le \quad CR^{d-p}\left(\int_{A_R(y)} \frac{(\delta_0+|\nabla u|^2)^{\frac{p-2}{2}} |\nabla u \cdot (x-y)|^2}{|x-y|^{d-p+2}} \dx \right)^{\frac12}  \left(  \int_{B_{2R}(y)} \frac{(\delta_0+|\nabla u|^2)^{\frac{p}{2}}}{R^{d-p}}\dx \right)^{\frac{1}{2}}.
\end{aligned}
$$
Finally, for the last term on the right-hand side of \eqref{C2}, we use the H\"{o}lder and the Poincar\'{e} inequalities together with \eqref{basis6} with $\gamma:=d-p-\frac{p(1-\theta)(p-1)}{2(1 +\theta (p -1))}$ and also the a~priori bound \eqref{E2} to deduce that
\begin{equation*}
\begin{aligned}
R^{-1}\int_{A_R(y)} (\delta_0+ |\nabla u|)^{\theta(p-1)} |u-\overline{u}_R| \dx &\le C R^{\frac{d(p-1)(1-\theta)}{p}}\left(\int_{B_{2R}(y)} (\delta_0 +|\nabla u|^2)^{\frac{p}{2}}\dx \right)^{\theta + \frac{1-\theta}{p}}\\
&\le C R^{d-p + \frac{(1-\theta)(p-1)}{2}} \left(\int_{B_{2R}(y)} \frac{(\delta_0 +|\nabla u|^2)^{\frac{p}{2}}}{R^{d-p-\frac{p(1-\theta)(p-1)} {2(1 +\theta (p -1))}}}\dx \right)^{\theta + \frac{1-\theta}{p}}\\
&\le CR^{d-p + \frac{(1-\theta)(p-1)}{2}} (1+\|f\|_2 + \|\partial_t u\|_2)^{p'}.
\end{aligned}
\end{equation*}

Substituting all the above estimates into \eqref{C2} and dividing the resulting inequality by $R^{d-p}$, we deduce
\begin{equation}\label{C3}
\begin{split}
\int_{B_{R}(y)} \frac{|\nabla u|^p}{R^{d-p}} \dx &\le  CR^{\frac{2p-d}{2}}(1+\|f\|_2 + \|\partial_t u\|_2) \left(\int_{B_{2R}(y)}\frac{|\nabla u|^p}{R^{d-p}} \dx \right)^{\frac{1}{p}}\\
&\quad+ C\left(\int_{A_R(y)} \frac{(\delta_0+|\nabla u|^2)^{\frac{p-2}{2}} |\nabla u \cdot (x-y)|^2}{|x-y|^{d-p+2}} \dx \right)^{\frac12}  \left(  \int_{B_{2R}(y)} \frac{(\delta_0+|\nabla u|^2)^{\frac{p}{2}}}{R^{d-p}}\dx \right)^{\frac{1}{2}}\\
&\quad + CR^{\frac{(1-\theta)(p-1)}{2}} (1+\|f\|_2 + \|\partial_t u\|_2)^{p'}\\
&\le \frac14  \int_{B_{2R}(y)}\frac{|\nabla u|^p}{(2R)^{d-p}} \dx + C\int_{A_R(y)} \frac{(\delta_0+|\nabla u|^2)^{\frac{p-2}{2}} |\nabla u \cdot (x-y)|^2}{|x-y|^{d-p+2}} \dx \\
&\quad +CR^{\alpha_2} (1+\|f\|_2 + \|\partial_t u\|_2)^{p'},
\end{split}
\end{equation}
where we used the Young inequality and defined
\begin{equation}
\alpha_2:= \min \left\{\frac{(2p-d)p}{2(p-1)}, \frac{(1-\theta)(p-1)}{2} \right\}.
\end{equation}

\subsection{Beyond the critical Morrey exponent} 
Here we combine the estimates from two previous sections and use the hole filling method. We divide the inequality \eqref{B2} by $C 2^{d-p+2}$ and add the result to the inequality \eqref{C3} and observe that for some $\varepsilon>0$ we have
\begin{equation} \label{D2}
\begin{split}
&\varepsilon\int_{B_R(y)} \frac{(\delta_0+|\nabla u|^2)^{\frac{p-2}{2}} |\nabla u \cdot (x-y)|^2}{ |x-y|^{d-p+2}} +\int_{B_{R}(y)} \frac{|\nabla u|^p}{R^{d-p}}\dx \\
&\qquad \le \frac12\int_{B_{2R}(y)} \frac{|\nabla u|^p}{(2R)^{d-p}} \dx+C\int_{A_R(y)} \frac{(\delta_0+|\nabla u|^2)^{\frac{p-2}{2}} |\nabla u \cdot (x-y)|^2}{|x-y|^{d-p+2}} \dx\\
&\qquad \ \qquad + CR^{\alpha_3}(1+\|f\|_2 + \|\partial_t u\|_2)^{p'},
\end{split}
\end{equation}
where $\alpha_3:=\min\{\alpha_1, \alpha_2\}$. Adding to both sides the corresponding integral in order to fill the hole on the right-hand side, we deduce that
\begin{equation} \label{D1}
\begin{split}
&(C+\varepsilon)\int_{B_R(y)} \frac{(\delta_0+|\nabla u|^2)^{\frac{p-2}{2}} |\nabla u \cdot (x-y)|^2}{ |x-y|^{d-p+2}} +\int_{B_{R}(y)} \frac{|\nabla u|^p}{R^{d-p}}\dx \\
&\qquad \le \frac12\int_{B_{2R}(y)} \frac{|\nabla u|^p}{(2R)^{d-p}} \dx+C\int_{B_{2R}(y)} \frac{(\delta_0+|\nabla u|^2)^{\frac{p-2}{2}} |\nabla u \cdot (x-y)|^2}{|x-y|^{d-p+2}} \dx\\
&\qquad \ \qquad  CR^{\alpha_3}(1+\|f\|_2 + \|\partial_t u\|_2)^{p'}.
\end{split}
\end{equation}
Denoting
$$
I_R:=(C+\varepsilon)\int_{B_R(y)} \frac{(\delta_0+|\nabla u|^2)^{\frac{p-2}{2}} |\nabla u \cdot (x-y)|^2}{ |x-y|^{d-p+2}} +\int_{B_{R}(y)} \frac{|\nabla u|^p}{R^{d-p}}\dx, 
$$
we see that \eqref{D1} implies that 
$$
I_R \le  CR^{\alpha_3}(1+\|f\|_2 + \|\partial_t u\|_2)^{p'} + \kappa I_{2R}, \qquad \kappa:= \min \left\{\frac12, \frac{C}{C+\varepsilon} \right\}<1.
$$
Itterating the above inequality, we deduce that for  
$$
\alpha_{\min}:= \frac{1}{2p} \min\{-\ln \kappa /\ln 2, \alpha_3 \}
$$
we have  
$$
I_R \le  C\left(\frac{R}{R_0}\right)^{p\alpha_{\min}} \left((1+\|f\|_2 + \|\partial_t u\|_2)^{p'} + I_{R_0} \right).
$$
Consequently, using also the definition of the set $K$ and the definition of $R_0$, as well as the a~priori estimate~\eqref{E2}, we deduce from the above inequality that
\begin{equation}
\label{kilo0}
\sup_{y\in K} \sup_{R\in (0,R_0)} \int_{B_R(y)} \frac{|\nabla u|^p}{R^{d-p+p\alpha_{\min}}} \dx \le C(K,R_0) \left(1+\|f\|_2 + \|\partial_t u\|_2\right)^{p'}.
\end{equation}
Combined with the estimate \eqref{E2} and the Morrey embedding, we obtain
\begin{equation}
\label{kilo}
\|u\|_{\mathcal{C}^{\alpha_{\min}}(K;\mathbb{R}^N)}  \le C(K,R_0) \left(1+\|f\|_2 + \|\partial_t u\|_2\right)^{\frac{1}{p-1}}.
\end{equation}

\subsubsection{No ellipticity condition}
In case we do not assume \eqref{5}, we have only \eqref{E2} at our disposal, and then it follows from \eqref{kilo0} and \eqref{kilo}, after integration over the time interval $(0,T)$, that
\begin{equation}\label{result1}
\int_0^T \|u(t)\|_{\mathcal{C}^{\alpha_{\min}}(K;\mathbb{R}^N)}^{2(p-1)} \dt \le C(\Omega,K,f,u_0)
\end{equation}
This is the claim of the main theorem when the ellipticity condition is not used.

\subsubsection{Ellipticity condition}
In case we do assume \eqref{5}, we can use \eqref{E2} and \eqref{E3} to deduce from \eqref{kilo0} and \eqref{kilo} that
\begin{equation}\label{result2}
\sup_{t\in (0,T)} \sup_{y\in K} \sup_{R\in (0,R_0)} \int_{B_{R(y)}}\frac{|\nabla u|^p}{R^{d-p+p\alpha_{\min}}}\dx \le C(\Omega,K,f,u_0)
\end{equation}
and also
\begin{equation}\label{result2.5}
\|u\|_{L^{\infty}(0,T; \mathcal{C}^{\alpha_{\min}}(K;\mathbb{R}^N))} \le C(\Omega,K,f,u_0).
\end{equation}
The estimate \eqref{result2} is the desired estimate of the main theorem. Moreover, thanks to \eqref{E3}, we have that $\partial_t u \in L^{\infty}(0,T; L^2(\Omega; \mathbb{R}^N))$, so we can use the cross interpolation. Then it follows from \eqref{result2.5} that there exists $\beta>0$, depending on the dimension $d$ and on $\alpha_{\min}$, such that
\begin{equation}
\label{result3}
\|u\|_{\mathcal{C}^{\beta}((0,T)\times K;\mathbb{R}^N)} \le C(\Omega,K,f,u_0),
\end{equation}
which finishes the proof of the main theorem.
%


\end{document}